\documentclass[12pt]{article}
\usepackage{amssymb,latexsym,amsmath}
\def\C{\mathbb C}
\def\R{\mathbb R}

\def\Z{\mathbb Z}

\def\T{\mathbb T}

\newtheorem{thm}{Theorem}[section]
\newtheorem{lem}{Lemma}[section]

\newtheorem{prop}{Proposition}[section]

\begin{document}
\sffamily 
\sffamily
\title{Complex flows, escape to infinity\\  and a question of Rubel }
\author{J.K. Langley}
\maketitle

\begin{abstract}
Let $f$ be a transcendental entire function. It was shown in a previous paper \cite{Latraj} 
that the holomorphic flow $\dot z = f(z)$ always has infinitely many trajectories tending to infinity in finite time. It will be proved here that such trajectories are in a certain sense rare,
although an example will be given to show that there can be uncountably many.
In contrast, for  the classical antiholomorphic flow 
$\dot z = \bar f(z)$, such trajectories need not exist at all, although they must if $f$
belongs to  
the Eremenko-Lyubich class $\mathcal{B}$.
It is also shown that for transcendental entire $f$ in $\mathcal{B}$
 there exists a path tending to infinity
on which $f$ and all its derivatives tend to infinity, thus affirming a conjecture of Rubel for this class.
\end{abstract}

\section{Introduction}

The starting point of this note is   the  flow 
\begin{equation}
 \label{H1}
\dot z = f(z),
\end{equation}
in which $f $ or its conjugate
 $\bar f$ is an entire function.  A trajectory for (\ref{H1})
is a path $z(t)$ in the plane with  $z'(t) = f(z(t)) \in \C$ for $t$ in some maximal interval $(\alpha, \beta)  \subseteq \R$. By the existence-uniqueness theorem, such trajectories are either constant (with $z(t)$ a zero of $f$), periodic or
injective.
It was shown in 
\cite[Theorem 5]{kingneedham}  that if  $f$ is a polynomial in $z$
of degree  $n \geq 2$  then there exist $n-1$ disjoint trajectories  for (\ref{H1}) which
 tend to
infinity in finite increasing time, that is, which satisfy
$\beta \in \R$ and $\lim_{t \to \beta - } z(t) = \infty$.The following theorem for holomorphic flows with transcendental entire $f$ was proved in \cite[Theorem 1.1]{Latraj}. 

\begin{thm}[\cite{Latraj}]
 \label{thm0}
Let the function $f$ be  transcendental entire: then (\ref{H1}) has infinitely many pairwise disjoint trajectories which 
tend to
infinity in finite increasing time.
\end{thm}

For meromorphic functions in general, such trajectories need not exist at all  \cite{Latraj},
but a result was also proved in \cite{Latraj} for the case where  
$f$ is  transcendental and meromorphic in the plane
and the inverse function
$f^{-1}$ has a logarithmic singularity over $\infty$: this means that there exist $M > 0$ and a component $U$ of the set $\{ z \in \C : \, |f(z)| > M \}$ such that $U$ contains no poles of $f$ 
and $\log f$ maps $U$ conformally onto the half-plane 
$H = \{ v \in \C : \, {\rm  Re } \, v > \log M \}$ \cite{BE,Nev}. 
 In this case \cite[Theorem 1.2]{Latraj},
(\ref{H1})  has infinitely many pairwise disjoint trajectories tending to 
infinity in finite increasing time from within a neighbourhood $\{ z \in U : |f(z)| > M' \geq M \}$ 
of the singularity.

On the other hand, for entire $f$ in (\ref{H1}),  it seems that  trajectories which  tend to
infinity in finite increasing time are somewhat exceptional. For 
the simple example $\dot z = - \exp( -z)$, it is easy to check that all trajectories satisfy
$\exp(z(t)) = \exp(z(0)) -t$ and so tend to infinity as $t$ increases, but take infinite time to do so
unless $\exp( z(0))$ is real and positive.  

It will be shown that for transcendental
entire $f$ there is, in a certain sense, zero probability of landing on a trajectory of (\ref{H1}) which tends to infinity in finite time.
To state the theorem, let  $f$ be transcendental
 entire and let
\begin{equation}
 \label{Fdef1}
z_0 \in \C, \quad f(z_0) \neq 0, \quad 
F(z) = \int_{z_0}^z \frac{du}{f(u)} .
\end{equation}
Then $F(z)$ is defined near $z_0$ and is real and increasing as $z$ follows the trajectory $\zeta_{z_0} (t)$ of (\ref{H1})
starting at $z_0$. Let $\delta $ be small and positive and take the pre-image $ L_\delta(z_0)$ of the real interval
$(- \delta, \delta)$ under the function $- i F(z) $; then $ L_\delta(z_0)$ is perpendicular to $\zeta_{z_0} (t)$ at $z_0$.
The proof of the following result is adapted from that of the Gross star theorem \cite[p.292]{Nev}. 

\begin{thm}
 \label{thmhol}
Let $f$ be a transcendental entire function and let $z_0$ and  $F$ be as in (\ref{Fdef1}). 
For small positive $\delta$ let $Y_\delta$ be the set of $y \in (- \delta, \delta)$ 
such that the trajectory of (\ref{H1}) starting at $F^{-1}(iy)$ tends to infinity in finite increasing time.
Then $Y_\delta$ has Lebesgue measure $0$. 
\end{thm}
Theorem \ref{thmhol}  seems unlikely to be best possible, but an example from \cite{Volk} (see \S \ref{uncountable}) shows that there exists a transcendental
 entire $f$ for which
 (\ref{H1}) has uncountably many trajectories tending  to infinity in finite increasing time.


It seems natural to ask similar questions in respect of the antiholomorphic flow 
\begin{equation}
 \label{AH}
\dot z = \frac{dz}{dt} =  \bar g(z),
\end{equation}
where $g$ is a non-constant entire function.
Equation (\ref{AH}) appears widely in textbooks as a model for incompressible irrotational plane fluid flow,
and is linked to (\ref{H1}) insofar as if $f = 1/g$ then (\ref{AH}) has the same trajectories as 
(\ref{H1}),  since $\bar g = f/|f|^2$, 
although zeros of one of $f$ and $g$
are  of course poles of the other and in general the speeds of travel 
differ.
The trajectories  of (\ref{AH})  are determined by choosing $G$ with 
$G'(z) = g(z)$ and writing
\begin{equation}
 \label{transform1}
v = G(z), \quad \dot v = g(z) \dot z = |g(z)|^2  \geq 0 ,
\end{equation}
which leads to the classical fact  that trajectories for (\ref{AH})
are level curves of ${\rm Im} \, G(z)$ on which ${\rm Re} \, G(z)$ increases with $t$. 
By the maximum principle,  ${\rm Im} \, G(z)$ cannot be constant on a closed curve.
Thus, apart from the countably many which tend to a zero of $G' = g$, 
all trajectories for (\ref{AH}) go to infinity, but this leaves open the  question as to how long  they take to do so.

If a non-constant trajectory $\Gamma$
of (\ref{AH}) passes from  $z_1 $ to $ z_2$ along an arc meeting no zeros of $g$, 
then ${\rm Im} \, v = \beta$ is
constant on $\Gamma$ and $X = {\rm Re} \,  v$ increases from $X_1
=  {\rm Re} \, G(z_1) $ to $X_2  =  {\rm Re} \, G(z_2)$. 
Thus (\ref{transform1}) implies that 
the transit time is
\begin{equation}
\int_{X_1+i\beta}^{X_2+i \beta } \frac1{|g(z)|^2}  \, dv =
 \int_{X_1+i \beta}^{X_2+i \beta} \left| \frac{dz}{dv} \right|^2  \, dv 
  =
 \int_{X_1}^{X_2} \left| \frac{dz}{dX} \right|^2  \, dX . 
\label{transit}
\end{equation}
This formula shows that a zero of $g$ cannot be reached in finite time, because if $z$
tends to a zero $z_3$ of $g$ of multiplicity $m$ as $X \to X_3$ then, with $c_j$ denoting non-zero constants,
\begin{eqnarray*}
X - X_3 &=& G(z)-G(z_3) 
\sim c_1 (z-z_3)^{m+1}, \\
\left| \frac{dz}{dX} \right|^2 &=&
\frac1{ |g(z)|^2 } \sim \frac{c_2}{  |X-X_3|^{2m/(m+1)} }
\geq \frac{ c_2 }{|X - X_3|} .
\end{eqnarray*}
Suppose now that $G' = g$ is a  polynomial of degree $n \geq 1$ in (\ref{AH}),
(\ref{transform1}) and (\ref{transit}).
If $S \in \R$ and $R$ is sufficiently large and positive then each pre-image  under $v = G(z)$
of the half-line $v = r + iS, r \geq R,$ gives a trajectory of  (\ref{AH}) 
which tends to infinity,
on which (\ref{transform1}) delivers
$$\frac{dt}{dv} = \frac1{|g(z)|^2} \sim \frac{c_3 }{ |z|^{2n}} \sim \frac{c_4}{  |v|^{2n/(n+1)}} .$$
Hence (\ref{transit}) implies that 
the transit time to infinity is finite for $n \geq 2$ and infinite for $n=1$. Thus, if $g$ is a non-linear polynomial, (\ref{AH}) always has uncountably many  trajectories tending to infinity in finite increasing time, but 
 this need not be the case for transcendental entire $g$. 

\begin{thm}
\label{thmbbh}
There exists a transcendental entire function $g$ such that (\ref{AH}) has no trajectories tending to infinity in finite increasing time. 
\end{thm}

Theorem \ref{thmbbh} also marks a sharp contrast with 
Theorem~\ref{thm0}, and its proof
rests on  the following immediate consequence of a result of Barth, Brannan and Hayman
\cite[Theorem 2]{BBH}. 

\begin{thm}[\cite{BBH}]
\label{BBHthm}
There exists a transcendental  entire function $G$ such that any unbounded connected plane set
contains a sequence $(w_n)$  tending to infinity on which 
$U = {\rm Re} \, G$ satisfies $(-1)^n U(w_n)  \leq 
 |w_n|^{1/2 } $.
\end{thm} 
To establish Theorem \ref{BBHthm}, it is only necessary to take the plane harmonic function $v$
constructed in \cite[Theorem 2]{BBH}, with the choice of $\psi(r)$ given by
\cite[p.364]{BBH}. With $U = v$, and $V$ a harmonic conjugate of $U$,
elementary 
considerations show that the resulting entire function $G = U+iV$ cannot be a polynomial.

On the other hand, in the presence of a logarithmic singularity of the inverse function 
over infinity,  trajectories of (\ref{AH})  tending to infinity in finite increasing time exist in
abundance.

\begin{thm}
 \label{thm2}
Let $g$ and $G$ be transcendental meromorphic functions in the plane such that $G'=g$ and
either $G^{-1}$ or  $g^{-1}$ has a logarithmic singularity over $\infty$. 
Then in each neighbourhood of the singularity   the flow (\ref{AH}) has a family of 
pairwise disjoint trajectories $\gamma_Y, Y \in \R$,  each of which tends to infinity in finite increasing time.
\end{thm}

Theorem \ref{thm2} applies in particular if $g$ or  its antiderivative 
$G$ is a transcendental entire function and 
belongs to the Eremenko-Lyubich class
$\mathcal{B}$,   which  plays a salient role in complex dynamics  \cite{Ber4,EL,sixsmithEL}
and 
is defined by the property that $F \in \mathcal{B}$
if the finite critical and asymptotic values of $F$ form a bounded set,
from which it follows that if $F \in \mathcal{B}$ is transcendental entire then  $F^{-1}$ automatically has  a logarithmic 
singularity over $\infty$.
A specific function to which 
Theorem \ref{thm2} may be applied is   $g(z) = e^{-z} + 1$; here $g$ is in $\mathcal{B}$, 
but its antiderivative $G$ is not,
and this example
also gives uncountably many trajectories  of   (\ref{AH}) taking infinite time to reach infinity
through the right half-plane. 

Theorem \ref{thm2} is quite straightforward to prove when the inverse of $G$ 
has a logarithmic singularity over infinity, but the method turns out to have a bearing
on the following question of Rubel
\cite[pp.595-6]{Linear}: if $f$ is a transcendental entire function, must there exist a path tending to
infinity on which $f$ and its derivative $f'$ both have asymptotic value $\infty$? This problem
 was motivated by the classical theorem of Iversen \cite{Nev},
which states that 
$\infty$ is an asymptotic value of every  non-constant entire function. 
For transcendental entire $f$ of finite order, a strongly affirmative answer to Rubel's question was provided by the 
following result \cite[Theorem 1.5]{Larubel}.

\begin{thm}[\cite{Larubel}]
\label{rubelthm}
Let the function $f$ be transcendental and meromorphic in the plane, of finite order of growth, 
and with finitely many poles.
Then there exists a path $\gamma$ tending to infinity such that,  for each non-negative integer $m$
and each positive real number $c$, 
\begin{equation}
\lim_{z \to \infty, z \in \gamma
 } \frac{ \log |f^{(m)}(z)|}{ \log |z|} = +  \infty \quad \hbox{and} \quad \int_\gamma |f^{(m)}(z)|^{-c} |dz| < + \infty .
\label{rr3}
\end{equation}
\end{thm}

For  functions of infinite order, Rubel's question appears to be difficult, 
although a path satisfying (\ref{rr3}) for $m=0$ is known to exist for any transcendental entire
function $f$ \cite{LRW}. 
However, a direct analogue of Theorem \ref{rubelthm} goes through  relatively straightforwardly
for  transcendental entire functions $f$  
in the Eremenko-Lyubich class $\mathcal{B}$.

\begin{thm}
 \label{thm1}
Let $f$ be a transcendental meromorphic function in the plane such that $f^{-1}$ has a logarithmic singularity over $\infty$, and let $D \in \R$. 
Then there exists a path $\gamma$ tending to infinity in a neighbourhood of the singularity, 
such that $f(z) -iD$ is real, positive and increasing on $\gamma$ and
(\ref{rr3})  holds for each integer $m \geq 0 $
and  real  $c > 0$. 
\end{thm}

 This paper is organised as follows: Theorem \ref{thmhol} is proved in 
 \S\ref{pfthmhol}, followed by an example in \S\ref{uncountable} and the proof
 of Theorem \ref{thmbbh} in \S\ref{pfthmbbh}. It is then convenient to give the proof
 of Theorem \ref{thm1} in \S\ref{pfthm1}, prior to that of Theorem \ref{thm2} in
 \S\ref{pfthm2}. 

\section{Proof of Theorem \ref{thmhol}}\label{pfthmhol}

Let $f$, $F$, $z_0$ and $\delta$  be as in the statement of Theorem \ref{thmhol}.
For  $y \in (- \delta, \delta)$ let  $g(y) = F^{-1}(iy)$ and let 
$T(y)$ be the supremum of $s > 0$ such that  the trajectory  $\zeta_{g(y)}(t)$ of (\ref{H1}) with $\zeta_{g(y)}(0) = g(y)$ is defined  and injective
for 
$0 \leq t < s$. 
If the trajectory $\zeta_{g(y)}(t)$ is periodic with minimal period $S_y$  then $T(y) = S_y$ 
and
$\zeta_{g(y')}(t)$ has the same period for $y'$ close to $y$
\cite{brickman}. 
Furthermore, if  $\zeta_{g(y)}(t)$ tends  to infinity in finite time then $T(y) < +  \infty$, while if
$T(y)$ is finite but  $\zeta_{g(y)}(t)$ is not periodic then 
$\lim_{t \uparrow T(y)} \zeta_{g(y)}(t) = \infty$ 
\cite[Lemma 2.1]{Latraj}. Set
$$
A =  \{ iy + t: \,  \, y \in (- \delta, \delta), \, 0 < t < T(y) \} , \quad B = \{ \zeta_{g(y)} (t) : \, y \in (- \delta, \delta), \, 0 < t < T(y) \}.
$$
Then $G( iy + t ) = \zeta_{g(y)} (t)   $
is a bijection from $A$ to $B$.

For  $u = \zeta_{g(y)} (t)$, where $y \in (- \delta, \delta)$ and $0 < t < T(y)$, let $\sigma_u$ be 
the subarc of 
$ L_\delta(z_0)$ from $z_0$ to $g(y)$ followed by the sub-trajectory of (\ref{H1})  from $g(y)$ to $u$, and 
define $F$ by (\ref{Fdef1}) on a simply connected neighbourhood $D_u$ 
of $\sigma_u$. Then $F$ maps $\sigma_u$ bijectively to the line segment $[0, iy]$ followed by 
the line segment $[iy, iy+t]$, and taking a sub-domain if necessary makes it
possible to assume that $F$ is univalent on $D_u$, with inverse function defined on a neighbourhood of $[iy, iy+t]$. 

Let  $y'$ and $t'$ be real and close to $y$ and $t$ respectively. Then
the image under $F^{-1}$ of the line segment $[iy', iy' + t']$ is
an injective sub-trajectory of (\ref{H1}) joining $g(y') \in  L_\delta(z_0)$  to $F^{-1}(iy' +t')
= \zeta_{g(y')} (t') = G(iy'+t')$, and so $T(y') \geq t'$. Thus
$y \rightarrow T(y)$ is lower semi-continuous and $A$ is a domain,
while $G: A \to B$ is analytic. Moreover, $A$ is
simply connected, because its complement in $\C \cup \{ \infty \}$
 is connected, and so is $B$. Furthermore,
 $F$ extends to be analytic on $B$, by (\ref{Fdef1}) and the fact
 that  $f \neq 0$ on $B$,
 and $F \circ G$ is the identity on $A$ because $F(G(t)) = t$ for small positive $t$. 

For $N \in (0,  + \infty ) $, let  $M_N $ be the set of all $  y$ in $ (- \delta, \delta)$ such that 
$\zeta_{g(y)} (t)$ tends to infinity and
$T(y) < N $. 
To prove Theorem \ref{thmhol}, it suffices to show that each such $M_N$ has measure $0$, and the subsequent steps will be adapted from the proof of the Gross star theorem \cite[p.292]{Nev} and its extensions due to Kaplan \cite{Kaplan}.
Let
$\Lambda_N \subseteq B$ be the image of $\Omega_N = \{ w \in A : \, {\rm Re} \, w < N \}$ under $G$, let $r$ be large and positive
and denote the circle $|z| = r$ by $S(0, r)$. 
Then $S(0, r) \cap \Lambda_N$ is a  union of countably many open arcs $\Sigma_r$. 

If $y \in M_N$ then $T(y) < N$ and 
as $t \to T(y)$ the image $z = G(iy+t) $ tends to infinity in $\Lambda_N$ and so crosses $S(0, r)$, 
and hence there exists $\zeta $ in some $ \Sigma_r$ with ${\rm Im} \, F(\zeta) = y$, since
$F: B \to A$ is the inverse of $G$. 
Thus the measure $\mu_N$ of $M_N$ is at most 
the total length $s(r)$ of the arcs $F(\Sigma_r)$. It follows  from the Cauchy-Schwarz inequality that, as $t \to + \infty$,  
\begin{eqnarray*}
  \mu_N^2 &\leq&   s(t)^2 =  \left( \int_{t e^{i \phi } \in \Lambda_N } |F'(t e^{i \phi } )| t \, d \phi \,    \right)^2 \\
&\leq&   \left( \int_{t e^{i \phi } \in \Lambda_N } |F'(t e^{i \phi } )|^2 t \, d \phi \,    \right)
\left( \int_{t e^{i \phi } \in \Lambda_N }  t \, d \phi \,    \right) \leq  2 \pi  t \left( \int_{t e^{i \phi } \in \Lambda_N } |F'(t e^{i \phi } )|^2 t \, d \phi \,    \right) .
\end{eqnarray*}
Thus  $\mu_N = 0$, since dividing by $2 \pi t$ and integrating from $r$ to $r^2$ yields,  as $r \to + \infty$,  
\begin{eqnarray*}
\frac{ \mu_N^2 \log r }{2 \pi} 
&\leq&    \int_r^{r^2}   \int_{t e^{i \phi } \in \Lambda_N } |F'(t e^{i \phi } )|^2 \, t \, d \phi \,   dt 
\leq  \int_{\Lambda_N}   |F'(t e^{i \phi } )|^2 \, t  \, d \phi \,  dt =      \hbox{area $(\Omega_N)$} \leq 2 \delta N .
\end{eqnarray*}
\hfill$\Box$
\vspace{.1in}

\section{An example}\label{uncountable}

Suppose that  $G$ is  a locally univalent meromorphic function  in the plane, 
whose set of asymptotic values  is an uncountable subset $E$ of  the unit circle $\T$. Suppose further that there exists 
a simply connected plane domain $D$, mapped univalently onto the unit disc $\Delta$ by $G$, 
such that the branch  $\phi$  of $G^{-1}$ mapping $\Delta$  to $D$ has no analytic
extension to a neighbourhood of any $\beta \in E$. 

Let $F = S(G)$, where $S$ is  a M\"obius transformation mapping $\Delta$ onto 
$\{ w \in \C : \, {\rm Re} \, w < 0 \}$, and for $\beta \in E$
let $\alpha = S(\beta)$ and let $L$ be the half-open line segment $[\alpha -1, \alpha)$. Then $M = S^{-1}(L)$ is 
a line segment or circular arc in $\Delta$ which meets $\T$ orthogonally at $\beta$. Moreover, $\phi(M)$ is a level curve of ${\rm Im} \, F$ in $D$, 
which cannot tend to a simple $\beta$-point of $G$ in $\C$ because this would imply that $\phi$ extends
to a neighbourhood of $\beta$. Hence $\phi(M)$ is a path tending to infinity in $D$, on which ${\rm Im} \, F(z)$ is constant
and $F(z)$ tends to $\alpha$.

Since $G$ and $F$ are locally univalent,  $f = 1/F'$ is entire. As $t \to 0-$ write, on $\phi(M)$,
$$
F(z) = \alpha + t, \quad \quad \frac{dt}{dz} = F'(z) = \frac1{f(z)}, \quad \frac{dz}{dt} = f(z),
$$
so that $\phi(M)$ is a trajectory of (\ref{H1}) which tends to infinity in finite increasing time, and there exists one of these for
every $\beta$ in the uncountable set $E$.

A suitable $G$  is furnished by a construction of Volkovyskii  \cite{Ermich,Volk}, 
in which $\T \setminus E$ is a union of  disjoint
open circular arcs $I_k = (a_k, b_k)$, oriented counter-clockwise.
For each $k$, take the  multi-sheeted
Riemann surface onto which $(a_k - b_k e^z)/(1-e^z)$ maps the plane, cut it along a curve  which projects to $I_k$,
and glue to $\Delta$  that half which lies to the right as $I_k$ is followed counter-clockwise. This forms a
simply connected Riemann surface $R$ with no algebraic branch points. 
By  \cite[Theorem 17, p.71]{Volk}
(see also \cite[p.6]{Ermich}),  the $I_k$ 
can be chosen so that $R$ is parabolic and is thereby the image surface of  a locally univalent meromorphic function $G$ in the plane.
\hfill$\Box$
\vspace{.1in}

\section{Proof of Theorem \ref{thmbbh}}\label{pfthmbbh} 

Following the notation of the introduction, suppose that $v=G(z)$ is a  transcendental entire function with derivative $g$ in (\ref{AH}), (\ref{transform1})
and  (\ref{transit}). 

\begin{prop}
\label{propbbh}
Let $\Gamma$ be a level curve tending to infinity 
on which $Y = {\rm Im} \, G(z) = \beta \in \R $ and
$X = {\rm Re} \, G(z) $ increases, with $X 
\geq \alpha \in \R $, and assume that $\Gamma$ meets no zero of $g$. 
Suppose that $(z_n)$ is a sequence tending to infinity on $\Gamma$ 
such that $v_n = G(z_n ) = X_n + i \beta $
satisfies $v_n = o( |z_n|)^2 $. 
Then 
the trajectory of (\ref{AH}) which follows $\Gamma$ takes infinite time in tending to infinity.
\end{prop}
Here it is not assumed or required 
 that $X \to + \infty$ as $z \to \infty$ on $\Gamma$. \\
 \\
 \textit{Proof of Proposition \ref{propbbh}.} 
 It may be assumed that $\Gamma$ starts at $z^*$ and $G(z^*) = \alpha + i \beta$. 
Denote positive constants, independent of $n$, by $C_j$. Then the Cauchy-Schwarz
inequality  gives,
 as $n $ and $z_n$ tend to infinity,
\begin{eqnarray*} 
|z_n|^2 &\leq&   \left( C_1 + \int_\alpha^{X_n} \left| \frac{dz}{dX} \right| \, dX \right)^2 \\
&\leq& 2 \left(  \int_\alpha^{X_n} \left| \frac{dz}{dX} \right| \, dX \right)^2 \\
&\leq&  2 \left(  \int_\alpha^{X_n}  \, dX \right)
 \left(  \int_\alpha^{X_n} \left| \frac{dz}{dX} \right|^2 \, dX \right) \\
 &\leq& 2  \left(  |v_n| + C_2 \right)
 \left(  \int_\alpha^{X_n} \left| \frac{dz}{dX} \right|^2 \, dX \right) \\
  &\leq& o  \left(  |z_n|^2 \right)
 \left(  \int_\alpha^{X_n} \left| \frac{dz}{dX} \right|^2 \, dX \right) .
 \end{eqnarray*}
Thus (\ref{transit}) shows that the transit time from
$z^*$ to $z_n$ tends  to infinity with $n$. 
\hfill$\Box$
\vspace{.1in}

\textit{Proof of Theorem \ref{thmbbh}.} 
Let $G$ be the  entire function  given by Theorem \ref{BBHthm}, and set $g = G'$. 
As noted in the introduction,
no trajectory of (\ref{AH}) can pass through a zero of $g$, and in any case it takes infinite time for a
trajectory  to approach a zero of $g$. Furthermore, if $\Gamma$ is 
a level curve, starting at $z^*$ say,
 on which ${\rm Im} \, G(z)$ is constant and $U(z) =  {\rm Re} \, G(z)$ 
increases, and on which $g$ has no zeros,  then there exists a  sequence $z_n = w_{2n}$ 
which tends to infinity on $\Gamma$ 
and  satisfies 
$$U(z^*) \leq U(z_n)  \leq |z_n|^{1/2} , \quad 
|G(z_n)| \leq |U(z_n)| + O(1) \leq   |z_n|^{1/2} + O(1).$$
Hence $\Gamma$ 
satisfies the hypotheses of Proposition \ref{propbbh}. 
It now follows that (\ref{AH}) has no trajectories tending to infinity in finite increasing time.
Since time can be reversed for these flows by setting $s = -t$ and 
$dz/ds = - \bar g(z)$, the same example has  no trajectories tending to infinity in finite
decreasing  time either. 
\hfill$\Box$
\vspace{.1in}


\section{Proof of Theorem \ref{thm1}}\label{pfthm1}

Let $f$ be as in the hypotheses. 
Then there exist $M > 0$ and a component $U$ of  $\{ z \in \C : \, |f(z)| > M \} $ such that $v = \log f(z)$ is a conformal bijection
from $U$ to the half-plane $H$
given by ${\rm Re} \, v >  N = \log M$; it may be assumed that $0 \not \in U$.
Let $\phi: H \to U$ be the inverse function. If $u \in H$ then $\phi$ and
$\log \phi$ are
 univalent on the disc $|w-u| < {\rm Re} \, u - N$ and so
Bieberbach's theorem and Koebe's quarter theorem \cite[Chapter 1]{Hay9}
imply that 
\begin{equation}
 \label{h3}
\left| \frac{\phi''(u)}{\phi'(u)} \right| \leq \frac4{{\rm Re} \, u -N } ,
\quad  \left| \frac{\phi'(u)}{\phi(u)} \right| \leq \frac{4 \pi}{{\rm Re} \, u -N }  .
\end{equation}

\begin{lem}
 \label{lem1} 
Let $v_0 $ be large and positive and 
for $0 \leq k \in \Z$ write 
\begin{equation}
 \label{rub1}
V_k = \left\{ v_0 + t e^{i \theta} : \, t \geq 0, \, - \,  \frac{\pi}{2^{k+2}} \leq \theta \leq \frac{\pi}{2^{k+2}} \right\}, \quad 
G_k(v) = \frac{f^{(k)}(z)}{f(z)}, \quad z = \phi(v). 
\end{equation}
Then there exist positive constants 
$d$ and  $c_k $ such that
$| \log \phi'(v) | \leq d  \log  ( {\rm Re} \, v )$ as $v \to \infty$ in $V_1$
and $| \log |G_k(v)| | \leq c_k \log  ( {\rm Re} \, v )$ as $v \to \infty$ in $V_k$. 
\end{lem}
\textit{Proof.} For $v \in V_1$,  parametrise the straight line segment from $v_0$ to $v$ with respect to $s = {\rm Re} \, u$. Then (\ref{h3}) and the simple estimate 
$|du| \leq \sqrt{2} ds $ yield  $| \log \phi'(v) | = O( \log  ( {\rm Re} \, v ))$ as $v \to \infty$ in $V_1$. 
Next, the assertion for $G_k$
 is trivially true for $k=0$, so assume that it holds for some $k \geq 0$ and write 
\begin{eqnarray*}
G_{k+1}(v) &=&
\frac{f^{(k+1)}(z)}{f(z)} = \frac{f^{(k)}(z)}{f(z)} \cdot \frac{f'(z)}{f(z)} +\frac{d}{dz} \left( \frac{f^{(k)}(z)}{f(z)} \right) \nonumber \\
&=& G_k(v) G_1(v) + \frac{G_k'(v)}{\phi'(v)} = \frac{G_k(v)}{\phi'(v)} \left( 1 + \frac{G_k'(v)}{G_k(v)} \right).  
\end{eqnarray*}
Thus it suffices to show that $G_k'(v)/G_k(v) \to 0$ as as $v \to \infty$ in $V_{k+1}$. By (\ref{rub1}) there exists a small positive
$d_1$ such that if $v \in V_{k+1}$ is large then the circle $|u - v| = r_v = d_1 {\rm Re} \, v$ lies in $V_k$, and  the
differentiated Poisson-Jensen formula \cite[p.22]{Hay2} delivers 
$$
\frac{G_k'(v)}{G_k(v)} = \frac1{ \pi} \int_0^{2 \pi} \, \frac{ \log | G_k(v + r_v e^{i \theta } )|}{r_v e^{i \theta}} \, d \theta 
= O \left( \frac{ \log  ( {\rm Re} \, v ) }{ {\rm Re} \, v } \right)  \to 0
$$
as $v \to \infty$ in $V_{k+1}$. This proves the lemma.
\hfill$\Box$
\vspace{.1in}

To establish Theorem \ref{thm1}, take any $D \in \R$. 
Then there exist 
$v_1 \in [1, + \infty) $ and a path 
$$\Gamma \subseteq 
\{ v \in \C : \, {\rm Re} \, v >  N , \, | {\rm Im} \, v | < \pi/4 \} \subseteq H$$
which is mapped by $e^v$ to the half-line 
$\{ t + iD : \, t \geq v_1 \}$. 
Thus $f(z) -i D = e^v -i D $ is real and positive for $z$ on $\gamma = \phi (\Gamma)$, 
and $\Gamma \setminus V_k $ is bounded for each $k \geq 0$. 
Now write, on $\Gamma$, 
$$
e^v = t+iD,
\quad \frac{dv}{dt} = \frac1{t+iD}, \quad 
s = {\rm Re} \, v 
= \frac12 \ln (t^2+D^2) .$$
Hence, for any non-negative integers $k, m$,  Lemma \ref{lem1} gives, 
as $v \to \infty$ on $\Gamma $,
$$
\left| \frac{f^{(k)}(z)}{z^m} \right| = \left| \frac{f(z) G_k(v)}{z^m} \right|   = \left| \frac{e^v G_k(v)}{\phi(v)^m} \right|
\geq \frac{e^s }{ s^{c_k+md}  } \geq e^{s/2}  \to \infty .
$$
It then follows that, for $c > 0$, 
\begin{eqnarray*}
\int_\gamma |f^{(k)}(z)|^{-c} \, |dz| &\leq & O(1) + 
\int_\Gamma e^{- cs/2} |\phi'(v)| \, |dv| \\
&\leq& O(1) + 
\int_\Gamma e^{- cs/4}  \, |dv| \\
&=& O(1) + 
\int_{v_1}^{+\infty} \frac1{(t^2+D^2)^{1/2+c/8}} \, dt  < + \infty  .
\end{eqnarray*}
\hfill$\Box$
\vspace{.1in}

\section{Proof of Theorem \ref{thm2}}\label{pfthm2}

Suppose first that the inverse function of 
the antiderivative $G$ of $g$ has a logarithmic singularity over infinity, and take $D \in \R$. 
Then  
Theorem \ref{thm1}
may be applied with $f = G$ and $m= c = 1$, giving a level curve  $\gamma = \gamma_D$, lying in a neighbourhood of the singularity,
on
which ${\rm Im} \, G(z) = D$ and  ${\rm Re} \, G(z) $ increases. This curve  is a trajectory for (\ref{AH}), 
traversed in time 
$$
\int_\gamma \frac1{\bar g(z)} 
\, dz  \leq \int_\gamma  |G'(z)|^{-1}  \, |dz| < + \infty ,$$
which completes the proof in this case. 

For the proof of the
 following lemma
 the reader is referred to the statement and proof of \cite[Lemma 3.1]{blnewqc}.

\begin{lem}[\cite{blnewqc}]
 \label{lemfirstest}
Let  the function $\phi : H \to \C \setminus \{ 0 \}$ 
be analytic and univalent, where $H = \{ v \in \C : \, {\rm Re} \, v > 0 \}$,
and
for $v, v_1 \in H$ define $Z(v) = Z(v, v_1)$  by
\begin{equation}
 \label{h1}
 Z(v, v_1) = \int_{v_1}^v e^{u/2} \phi'(u) \, du = 
2 e^{v/2} \phi'(v) - 2 e^{v_1/2} \phi'(v_1) - 2 \int_{v_1}^v e^{u/2} \phi''(u) \, du .
\end{equation}
Let 
$\varepsilon $  be a small positive real number. Then there exists a large positive real number $N_0$, depending  on $\varepsilon$ but not on $\phi$,
with the following property.

Let $v_0 \in H$ be such that  $S_0 = {\rm Re} \, v_0 \geq N_0$,
and define $v_1, v_2, v_3, K_2$ and $ K_3$ by 
\begin{equation*}
\label{vjdef}
v_j = \frac{2^j S_0}{128} + i T_0, 
 \quad T_0 = {\rm Im} \, v_0, \quad 
K_j = \left\{ v_j + r e^{i \theta} : \, r \geq 0, \, - \frac{\pi}{2^j} \leq \theta \leq  \frac{\pi}{2^j} \right\}.
\end{equation*}
Then the following two conclusions both hold:\\
(i)  $Z = Z(v, v_1)$ satisfies,  for $v \in K_2$,
\begin{equation}
 \label{h2}
 Z(v, v_1 ) = \int_{v_1}^v e^{u/2} \phi'(u) \, du  = 2 e^{v/2} \phi'(v) (1 + \delta (v) ), \quad | \delta (v) | < \varepsilon .
\end{equation}
(ii) $\psi = \psi (v, v_1) = \log Z(v, v_1) $
is univalent on a domain $H_1$,
with $v_0 \in H_1 \subseteq K_3$,  and $\psi(H_1)$ contains the strip 
\begin{equation}
 \label{Omegaimage}
\left\{  \psi (v_0)  + \sigma+ i \tau  : \, \sigma \geq \log \frac18 \, , \, -  2 \pi \leq  \tau \leq 2 \pi \right\} .
\end{equation}
\end{lem}
\hfill$\Box$
\vspace{.1in}

Assume henceforth that $g$ is as in the hypotheses of Theorem \ref{thm2}
and the inverse function of  $g$ has a logarithmic singularity over infinity. 
This time there exist $M > 0$ and a component $C$ of  
$\{ z \in \C : \, |g(z)| > M \} $ such that $\zeta  = \log g(z)$ is a conformal mapping of
$C$ onto the half-plane 
given by ${\rm Re} \, \zeta  > \log M$. Since 
(\ref{AH}) may be re-scaled via 
$z = Mw$ and $g(z) = Mh(w) $,
it may be assumed that $M = 1$ and $0 \not \in C$.
In order to apply Lemma \ref{lemfirstest},
let $\phi: H \to C$ be the inverse function $z = \phi(v)$ of the mapping
from $C$ onto $H$ given by
$$
v = 2 \zeta = 2 \log g(z), \quad g(z) = e^{v/2} ,
$$
As in the proof of Theorem \ref{thm1}, 
(\ref{h3}) holds for $u \in H$, with $N = 0$.

By (\ref{Omegaimage}) there exists $X_0 > 0$ such that 
$ Z(v, v_1)$ maps a domain $H_2 \subseteq H_1 \subseteq K_3 \subseteq H$
 univalently onto a half-plane ${\rm Re} \, Z > X_0$. Hence, for any $Y_0 \in \R$, 
there exists a path $\Gamma $ which tends to infinity in $ H_1 \subseteq K_3$ and
 is mapped by $ Z(v, v_1)$ onto 
the half-line $L_0 = \{  X + i Y_0, \, X \geq X_0 + 1 \}$.
Consider the flow
in $H_2$ given by 
\begin{equation}
 \label{vflow}
\phi'(v) \dot v = \overline{e^{v/2}} ;
\end{equation}
by (\ref{h2}) this transforms under $Z = Z(v, v_1)$ to
\begin{equation}
\label{wflow}
\dot Z = \frac{dZ}{dv} \, \dot v 
 = e^{v/2} \phi'(v) \dot v = | e^{v} | .
\end{equation}
Combining (\ref{h3}) and (\ref{h2}) shows that $ | e^{v} | \geq |Z(v)|^{3/2}   $ for large $v$ on $\Gamma$.  
Hence there exists a trajectory of (\ref{wflow}) which starts 
at $X_0+1 + iY_0$ and tends to infinity along $L_0$  in time 
$$
T_0 \leq  \int_{X_0+1}^\infty \left| \frac{dt}{dX} \right|  \, dX 
\leq O(1) + \int_{X_0+1}^\infty (X^2 + Y_0^2)^{-3/4} \, dX < + \infty .
$$
This gives a trajectory of (\ref{vflow}) tending to infinity along $\Gamma$ and taking finite time to do so, 
and hence a trajectory $\gamma$ of (\ref{AH}) in $C$, tending to infinity in finite increasing time. 
Since $Y_0 \in \R$ may be chosen at will, this proves Theorem \ref{thm2}. 

\hfill$\Box$

\noindent
School of Mathematical Sciences, University of Nottingham, NG7 2RD.\\
james.langley@nottingham.ac.uk

\end{document}